\documentclass{smfart}
\usepackage{bull,sabbah_dist-hol}

\begin{document}
\frontmatter
\title[D\'eveloppement asymptotique de~distributions~holonomes]
{D\'eveloppement asymptotique de~distributions~holonomes d'une~variable~complexe}

\alttitle{Asymptotic expansion of holonomic distributions of one complex variable}

\author[C.~Sabbah]{Claude Sabbah}
\address{UMR 7640 du CNRS\\
Centre de Math\'ematiques Laurent Schwartz\\
\'Ecole polytechnique\\
F--91128 Palaiseau cedex\\
France}
\email{sabbah@math.polytechnique.fr}
\urladdr{http://www.math.polytechnique.fr/~sabbah}

\begin{abstract}
Nous donnons la forme g\'en\'erale d'un germe de distribution holonome d'une variable complexe.
\end{abstract}

\begin{altabstract}
We give the general form for a germ of holonomic distribution of one complex variable.
\end{altabstract}


\altkeywords{Distribution, asymptotic expansion, holonomic, D-module}

\maketitle
\mainmatter

\section*{Introduction}
Dans cet article, nous utilisons la notion de \og dual hermitien d'un $\cD$-module\fg introduite par M\ptbl Kashiwara \cite{Kashiwara86} (voir aussi \cite{Bjork93}) pour donner la forme g\'en\'erale d'un germe de distribution d'une variable complexe satisfaisant \`a une \'equation diff\'erentielle holomorphe. Le cas o\`u l'\'equation est \`a singularit\'e r\'eguli\`ere est bien connu (voir par exemple \cite{B-M87}). Nous utilisons ici le fait que le dual hermitien d'un $\cD$-module holonome d'une variable complexe est encore holonome: c'est un cas particulier d'une conjecture g\'en\'erale de M.~Kashiwara; ce cas particulier est montr\'e dans \cite{Bibi97}, en analysant la dualit\'e hermitienne au niveau des cocycles de Stokes.

Nous n\'egligerons ci-dessous les distributions \`a support ponctuel (masses de Dirac) et travaillerons avec les germes de distributions mod\'er\'ees.

\section{D\'eveloppement asymptotique de distributions holonomes}\label{sec:asympt}
Soit $X$ un disque centr\'e en $0$ dans $\CC$, muni de la coordonn\'ee~$x$. Soit $u$ un germe en $0$ de distribution holonome mod\'er\'ee en $0$ sur $X$: autrement~dit,
\begin{enumerate}
\item
il existe un voisinage ouvert $U$ de $0$ dans $X$ tel que $u$ soit une distribution sur~$U^*=U\moins\{0\}$ qui soit la restriction d'une distribution sur $U$ (l'espace correspondant est not\'e $\Dbmo(U)$),
\item
il existe un op\'erateur diff\'erentiel lin\'eaire holomorphe non identiquement nul $P\in\cO(U)\langle\partial_{x}\rangle$ tel que l'on ait $P\cdot u=0$ dans $\Dbmo(U)$.
\end{enumerate}
On peut supposer que $U$ est choisi de sorte que $P$ n'ait de singularit\'e qu'en $0\in X$.

Soit $\pi:Y\to X$, $y\mto y^q=x$, un rev\^etement ramifi\'e de degr\'e $q\in\NN^*$. Alors l'image inverse par $\pi$ d'une distribution mod\'er\'ee en $0$ est bien d\'efinie comme distribution mod\'er\'ee en $0$ sur $Y$. Si $u$ est holonome, $\pi^*u$ l'est aussi.

\begin{theoreme}
Soit $u$ un germe en $0$ de distribution holonome mod\'er\'ee sur~$X$. Alors il existe:
\begin{itemize}
\item
un entier $q$, donnant lieu \`a un rev\^etement ramifi\'e $\pi:Y\to X$,
\item
un ensemble fini $\Phi\subset \ym\CC[\ym]$,
\item
pour tout $\varphi\in\Phi$, un ensemble fini $B_\varphi\in\CC$ et un entier $L_\varphi\in\NN$,
\item
pour tous $\varphi\in\Phi$, $\beta\in B_\varphi$ et $\ell=0,\dots,L_\varphi$, une fonction $f_{\varphi,\beta,\ell}\in \cC^\infty(Y)$
\end{itemize}
tels que l'on ait, dans $\Dbmo(V)$ et en particulier dans $C^\infty(V^*)$ (o\`u $V$ est un voisinage assez petit de $0$ dans $Y$), l'\'egalit\'e
\[\tag{$*$}
\pi^*u=\sum_{\varphi\in\Phi}\sum_{\beta\in B_\varphi}\sum_{\ell=0}^{L_\varphi}f_{\varphi,\beta,\ell}(y)e^{\varphi-\ov\varphi}\my^{2\beta}\Ly^\ell,
\]
o\`u on a not\'e
\[
\Ly\defin\module{\log\my^2}.
\]
\end{theoreme}

Notons que, pour $\varphi\in\ym\CC[\ym]$, la fonctions $e^{\varphi-\ov\varphi}$ est un multiplicateur dans $\Dbmo(V)$ (car c'est une fonction $C^\infty$ sur $V^*$, \`a croissance mod\'er\'ee \`a l'origine ainsi que toutes ses d\'eriv\'ees). Il en est de m\^eme des fonctions $\my^{2\beta}$ et $\Ly^\ell$.

Avant de montrer le th\'eor\`eme, nous allons pr\'eciser les $\varphi,\beta$ tels que $f_{\varphi,\beta,\ell}\neq0$ pour un certain~$\ell$, en rappelant d'abord des r\'esultats classiques sur la structure des connexions m\'eromorphes d'une variable (\cf \cite{Malgrange91} par exemple).

Soit $M$ un germe en $x=0$ de fibr\'e m\'eromorphe muni d'une connexion $\nabla$, \ie un $\CC\{x\}[\xm]$-espace vectoriel de dimension finie muni d'une connexion. Soit $\pi:y\mto x=y^q$ une ramification telle que le formalis\'e ramifi\'e $\wh N\defin\pi^+\wh M$ soit isomorphe au formalis\'e d'un fibr\'e m\'eromorphe \`a connexion \'el\'ementaire, \ie de la forme $N^{\el}=\oplus_{\varphi\in\Phi}(E^\varphi\otimes R_\varphi)$, o\`u les $R_\varphi$ sont \`a singularit\'e r\'eguli\`ere et~$E^\varphi$ est \'egal \`a $\cO_X$ muni de la connexion telle que $\nabla 1=d\varphi$; on note, pour plus de clart\'e, $\ephi$ la section $1$ de $E^\varphi$. Le groupe $\ZZ/q\ZZ$ agit naturellement sur~$\wh N$.

On note $V^\cbbullet M$ la filtration d\'ecroissante (index\'ee par $\RR$) de Deligne du fibr\'e m\'eromorphe $M$: chaque $V^bM$ est de type fini sur $\CC\{x\}\langle x\partial_x\}$, et le gradu\'e $\gr^b_VM\defin V^bM/V^{>b}M$ est un espace vectoriel de dimension finie sur lequel l'endomorphisme induit par $x\partial_x$ a ses valeurs propres~$\beta$ de partie r\'eelle \'egale \`a~$b$. On note alors $\psi_x^\beta M\subset\gr^b_VM$ le sous-espace propre g\'en\'eralis\'e associ\'e \`a~$\beta$. La multiplication par $x^k$ induit un isomorphisme $\gr^b_VM\isom\gr^{b+k}_VM$ qui transforme $x\partial_x$ et $x\partial_x+k$, donc qui induit un isomorphisme $\psi_x^\beta M\isom\psi_x^{\beta+k} M$.

La construction peut aussi \^etre appliqu\'ee au formalis\'e $\wh M$ et il est connu que $\gr^b_V\wh M=\gr^b_VM$ et $\psi_x^\beta\wh M=\psi_x^\beta M$.

Par ailleurs, si $N=\pi^+M$ comme ci-dessus, on a $V^bM=(V^{qb}N)^{\ZZ/q\ZZ}$ (on utilise le fait que, sur $1\otimes M$, $y\partial_y$ agit comme $q\id\otimes x\partial_x$). Enfin, en utilisant la d\'ecomposition de $\wh N=\wh N^\el$, on a $V^b(E^\varphi\otimes R_\varphi)=E^\varphi\otimes R_\varphi$ pour tout $b\in\RR$ si $\varphi\neq0$. Ainsi, $\psi^\beta_yN=\psi^\beta_y\wh N=\psi^\beta_yR_0$.

Revenons maintenant \`a la situation du th\'eor\`eme. Dans la suite, on travaille sur un voisinage assez petit de $0$ qu'on ne pr\'ecise pas, et qu'on note toujours $X$ ou $Y$.

Soit $M$ le $\cD_{X}[\xm]$-module engendr\'e par $u$ dans $\Dbmo(X)$. Alors $M$ est un $\cO_{X}[\xm]$-module localement libre de rang fini muni d'une connexion, induite par l'action de $\partial_x$. Soit $\pi:Y\to\nobreak X$ une ramification telle que $N\defin\pi^+M$ soit formellement isomorphe \`a $N^{\el}=\oplus_{\varphi\in\Phi}(E^\varphi\otimes\nobreak R_\varphi)$. On identifie le germe de fibr\'e \`a connexion $\pi^+N$ au $\cD_Y[\ym]$-sous-module de $\Dbmo(Y)$ engendr\'e par la distribution mod\'er\'ee $v=\pi^*u$. Par ailleurs, $N$ est de la forme $\cD_Y/(Q)$ pour un certain op\'erateur diff\'erentiel holomorphe $Q$ qui annule $v$.

\begin{definition}\label{def:sansram}
On dira que la distribution holonome mod\'er\'ee $u$ est \emph{sans ramification} si on peut choisir $\pi=\id$ ci-dessus.
\end{definition}

On travaillera directement avec la distribution holonome mod\'er\'ee sans ramification $v=\pi^*u$. Pour $\varphi\in\ym\CC[\ym]$, soit $N_\varphi$ le $\cD_Y[\ym]$-module engendr\'e par $e^{\ov\varphi-\varphi} v$ dans $\Dbmo{Y}$. C'est aussi un $\cO_{X}[\xm]$-module localement libre de rang fini muni d'une connexion.

Pour tout $\beta\in\CC$, on note $L'_{\varphi,\beta}(v)$ l'ordre de nilpotence de $y\partial_y-\beta$ sur $\psi_y^\beta(N_\varphi)$. On a bien s\^ur $L'_{\varphi,\beta}(v)=L'_{\varphi,\beta+k}(v)$ pour tout $k\in\ZZ$. Il existe alors un ensemble fini minimal $B'_\varphi(v)\subset\CC$ tel que, pour tout $j\in\NN$ on puisse trouver un entier $k(j)$ et un op\'erateur $P_j\in\CC\{y\}\langle y\partial_y\rangle$ tels que
\begin{equation}\label{eq:Bernsteinj}
\bigg[\prod_{k=0}^{k(j)}\prod_{\beta\in B'_\varphi(v)}\big[-(y\partial_y-\beta-k)\big]^{L'_{\varphi,\beta}}-y^jP_j\bigg]\cdot e^{\ov\varphi-\varphi} v=0.
\end{equation}

\begin{remarque}
Pour presque tout $\varphi$, on~a $B'_\varphi(v)=\emptyset$. On note $\Phi(v)$ l'ensemble des~$\varphi$ pour lesquels la composante $E^\varphi\otimes R'_\varphi$ de $\cD_Y[\ym]\cdot v$ n'est pas nulle. C'est aussi l'en\-semble des $\varphi$ pour lesquels la composante $\ov{E^{-\varphi}\otimes R''_\varphi}$ de $\cD_{\ov Y}[1/\ovy]\cdot v$ n'est pas nulle.
\end{remarque}

On d\'efinit de mani\`ere conjugu\'ee les objets $L''_{\varphi,\beta}$ et $B''_\varphi(v)$. On pose alors
\[
B_\varphi(v)=\Big[(B'_\varphi(v)-\NN)\cap B''_\varphi(v)\Big]\cup\Big[B'_\varphi(v)\cap(B''_\varphi(v)-\NN)\Big].
\]
Autrement dit, $\beta\in B_\varphi(v)$ si et seulement si $\beta\in B'_\varphi(v)\cup B''_\varphi(v)$ et $(\beta+\NN)\cap B'_\varphi(v)\neq\emptyset$ et $(\beta+\NN)\cap B''_\varphi(v)\neq\emptyset$. Pour tout $\beta\in\CC$, on pose $L_{\varphi,\beta}(v)=\min\{L'_{\varphi,\beta}(v),L''_{\varphi,\beta}(v)\}$, et on a $L_{\varphi,\beta+k}(v)=L_{\varphi,\beta}(v)$ pour tout $k\in\ZZ$.

Enfin, si $f\in C^\infty(Y)$, on d\'eveloppe $f$ par rapport \`a $y,\ovy$ et on peut associer \`a ce d\'eveloppement un ensemble minimal $E(f)\subset\NN^2$ tel que $f=\sum_{(\nu',\nu'')\in E(f)}y^{\nu'}\ovy^{\nu''}\!g_{(\nu',\nu'')}$ avec $g_{(\nu',\nu'')}\in C^\infty(Y)$. La minimalit\'e de $E(f)$ implique en particulier que $g_{(\nu',\nu'')}(0)\neq0$ si $(\nu',\nu'')\in E(f)$. On convient que $E(f)=\emptyset$ si $f$ est infiniment plate en $0$.

\begin{corollaire}\label{cor:asympt}
Soit $v$ une distribution holonome mod\'er\'ee sans ramification. Alors $v$ admet un d\'eveloppement $(*)$ dans $\Dbmo_Y$, avec $\Phi=\Phi(v)$ et $\beta\in B_\varphi(v)$. De plus, si $f_{\varphi,\beta,\ell}\neq0$ et si le point $(k',k'')\in\NN^2$ est dans $E(f_{\varphi,\beta,\ell})$, alors $\beta+k'\in B'_\varphi(v)+\NN$ et $\beta+k''\in B''_\varphi(v)+\NN$.
\end{corollaire}

\begin{proof}
Nous supposons le th\'eor\`eme d\'emontr\'e. On utilise la transformation de Mellin pour raisonner sur chaque coefficient du d\'eveloppement~$(*)$. Soit $\chi$ une fonction $C^\infty$ \`a support compact contenu dans un ouvert o\`u $v$ est d\'efinie, identiquement \'egale \`a $1$ pr\`es de $0$. On note de la m\^eme mani\`ere la forme $\chi\tfrac{i}{2\pi}\,dy\wedge d\ovy$. Par ailleurs, choisissons une distribution $\wt v$ induisant~$v$ sur~$Y^*$ et soit $p$ son ordre sur le support de $\chi$. Nous allons d'abord consid\'erer les coefficients pour lesquels $\varphi=0$.

Pour tous $k',k''\in\NN$, la fonction $s\mto\langle \wt v,\my^{2s}y^{-k'}\ovy^{-k''}\chi\rangle$ est d\'efinie et holomorphe sur le demi-plan $2\reel s>p+k'+k''$. Pour tout $j\geq1$, notons~$Q_j$ l'op\'erateur apparaissant dans \eqref{eq:Bernsteinj} (pour $\varphi=0$). Alors $Q_j\cdot \wt v$ est \`a support l'origine. Il sera commode dans la suite d'utiliser la notation $\alpha$ pour $-\beta-1$ et noter $A'_\varphi(v)=\{\alpha\mid\beta=-\alpha-1\in B'_\varphi(v)\}$. On en d\'eduit alors que, sur un demi-plan $\reel s\gg0$, la fonction
\[
\bigg[\prod_{k=0}^{k(j)}\prod_{\alpha\in A'_0(v)}(s-\alpha-k'+k)^{L'_{0,\alpha}}\bigg]\langle \wt v,\my^{2s}y^{-k'}\ovy^{-k''}\chi\rangle
\]
co\"incide avec une fonction holomorphe sur $2\reel s>p+k'+k''-j$. En appliquant le m\^eme argument de mani\`ere anti-holomorphe, on trouve que, pour tous $k',k''\in\NN$, la fonction $s\mto\langle \wt v,\my^{2s}y^{-k'}\ovy^{-k''}\chi\rangle$ s'\'etend en une fonction m\'eromorphe sur $\CC$ \`a p\^oles contenus dans $(A'_0(v)+k'-\NN)\cap (A''_0(v)+k''-\NN)$, l'ordre du p\^ole en $\alpha+\ZZ$ \'etant major\'e par $L_{0,\alpha}(v)$. De plus, cette fonction ne d\'epend pas du choix du rel\`evement $\wt v$.

Calculons maintenant la transform\'ee de Mellin du d\'eveloppement $(*)$ pour~$v$.

\begin{lemme}\label{lem:sanspole}
Si $\varphi\neq0$, alors, pour toute fonction $g\in\cC^\infty(Y)$, la transform\'ee de Mellin de $g(y)e^{\varphi-\ov\varphi}\my^{2\beta}\Ly^\ell$ est une fonction enti\`ere.
\end{lemme}

\begin{proof}
On montre que cette transform\'ee de Mellin est holomorphe sur tout demi-plan $\reel s>-p$ ($p\in\NN$). Pour cela, pour $p$ fix\'e, on d\'ecompose $g$ comme la somme d'un polyn\^ome en $y,\ov y$ et d'un reste qui s'annule \`a un ordre assez grand \`a l'origine pour que la partie correspondante de la transform\'ee de Mellin soit holomorphe sur $\reel s>-p$. On est donc ramen\'e \`a supposer que $g$ est un mon\^ome en $y,\ov y$. Il existe alors des \'equations fonctionnelles holomorphe et anti-holomorphe pour la distribution mod\'er\'ee $g(y)e^{\varphi-\ov\varphi}\my^{2\beta}\Ly^\ell$ avec polyn\^ome de Bernstein \'egal \`a $1$. On peut ainsi utiliser le m\^eme argument que ci-dessus, avec un terme entre crochets \'egal \`a $1$.
\end{proof}

Consid\'erons donc les termes du d\'eveloppement $(*)$ de $v$ pour lesquels $\varphi=0$. Il n'est pas restrictif de supposer que deux \'el\'ements distincts de l'ensemble d'indices $B_0$ qui intervient dans $(*)$ ne diff\`erent pas d'un entier, et que tout \'el\'ement $\beta$ de $B_0$ est maximal, en ce sens que l'escalier $\bigcup_\ell E(f_{0,\beta,\ell})$ est contenu dans $\NN^2$ et dans aucun $(m,m)+\NN^2$ avec $m\in\NN^*$. Soit $\beta\in B_0$. Nous utiliserons le fait que, pour tous $(\nu',\nu'')\in\ZZ^2$ \emph{non tous deux strictement n\'egatifs} et toute fonction $g\in C^\infty(Y)$ telle que $g(0)\neq0$, la fonction m\'eromorphe $s\mto\langle g(y)\my^{2\beta}\Ly^\ell,\my^{2s}y^{\nu'}\ovy^{\nu''}\chi\rangle$ a ses p\^oles contenus dans $\alpha-\NN$ (avec $\alpha=-\beta-1$), et a un p\^ole en $\alpha$ si et seulement si $\nu'=0$ et $\nu''=0$, ce p\^ole \'etant alors d'ordre $\ell+1$ exactement.

Pour $\beta\in B_0$, soit $E_\beta\subset\NN^2$ un ensemble minimal tel que $E_\beta+\NN^2=\bigcup_\ell (E(f_{0,\beta,\ell})+\NN^2)$. On d\'eduit de ce qui pr\'ec\`ede et du d\'eveloppement $(*)$ que, pour tout $(k',k'')\in E_\beta$, la fonction $s\mto\langle \wt v,\my^{2s}y^{-k'}\ovy^{-k''}\chi\rangle$ a un p\^ole non trivial en $\alpha$; de la premi\`ere partie de la preuve on conclut que $\alpha-k'\in A'_0(v)-\NN$ et $\alpha-k''\in A''_0(v)-\NN$, c'est-\`a-dire $\beta+k'\in B'_0(v)+\NN$ et $\beta+k''\in B''_0(v)+\NN$. Par hypoth\`ese sur $B_0$, il existe $(k',k'')\in E_\beta$ avec $k'=0$ ou $k''=0$. Il en r\'esulte que $\beta\in B_\varphi(v)$ et que la condition donn\'ee dans l'\'enonc\'e du corollaire est satisfaite par les \'el\'ements de $E_\beta$. Elle est alors aussi trivialement satisfaite par les \'el\'ements de tous les $E(f_{0,\beta,\ell})$.

Pour obtenir le r\'esultat pour les $f_{\varphi,\beta,\ell}$, on applique le r\'esultat pr\'ec\'edent \`a la distribution mod\'er\'ee $e^{\ov\varphi-\varphi}v$.
\end{proof}

\begin{proof}[\proofname\ du th\'eor\`eme]

Nous utiliserons le r\'esultat suivant:

\begin{theoreme}[{\cite[Prop\ptbl II.3.2.5]{Bibi97}}]
Pour $M$ comme ci-dessus, le $\cO_{\ov X}[\ov x^{-1}]$-module $\cHom_{\cD_X}(M,\Dbmo_X)$ est libre (de m\^eme rang que $M$) et muni d'une connexion anti-holomorphe, donc est un $\cD_{\ov X}$-module holonome.\qed
\end{theoreme}

On note $\Cmo_XM=\cHom_{\cD_X}(M,\Dbmo_X)$. On a donc un accouplement $\cD_X\otimes_\CC\cD_{\ov X}$-lin\'eaire canonique
\[
k:M\otimes_\CC\Cmo_XM\to\Dbmo_X,\qquad (m,\varphi)\mto\varphi(m).
\]
Puisque $M$ est engendr\'e par $u$ sur $\cD_X[\xm]$, $\varphi\in\Cmo_XM$ est d\'etermin\'e par sa valeur $\varphi(u)\in\Dbmo_X$. Il existe donc une section $\bun_u$ de $\Cmo_XM$ telle que $\bun_u(u)=u$.

Tout revient ainsi \`a montrer le th\'eor\`eme dans le cas o\`u
\[
k:M'\otimes_\CC\ov{M''}\to\Dbmo_X
\]
est un accouplement sesquilin\'eaire entre deux $\cO_X[\xm]$-modules libres de rang fini \`a connexion, $m',m''$ en sont deux sections locales, et $u=k(m',\ov{m''})$.

On se ram\`ene, par un rev\^etement cyclique, au cas o\`u $M'$ et $M''$ admettent chacun une d\'ecomposition formelle model\'ee sur $M^{\prime\el}$ et $M^{\prime\prime\el}$ (si $M''=\Cmo_XM'$, un rev\^etement qui convient pour l'un convient aussi pour l'autre).

On travaille ensuite avec l'image inverse par l'\'eclatement r\'eel $e:\wt Y\to Y$ de l'origine. On note $\cA_{\wt Y}$ le faisceau des fonctions $C^\infty$ sur $\wt Y$ annul\'ees par l'op\'erateur de Cauchy-Riemann $\ov{y\partial_y}$, $\cD_{\wt Y}=\cA_{\wt Y}\otimes_{e^{-1}\cO_Y}e^{-1}\cD_Y$ et $\Dbmo_{\wt Y}$ le faisceau sur $\wt Y$ des distributions mod\'er\'ees le long de $e^{-1}(0)=S^1$. On note enfin $\wt M=\cA_{\wt Y}\otimes_{e^{-1}\cO_Y}e^{-1}M$. C'est un $\cD_{\wt Y}$-module \`a gauche, qui est $\cA_{\wt Y}[\ym]$-libre.

L'accouplement $k$ s'\'etend de mani\`ere unique en un accouplement $\cD_{\wt Y}\otimes_\CC\cD_{\wt{\ov Y}}$-lin\'eaire
\[
\wt k:\wt M'\otimes_\CC\ov{\wt M''}\to\Dbmo_{\wt Y}
\]
(simplement parce que $M$ est $\CC\{y\}[\ym]$-libre). On alors peut travailler localement sur $\wt Y$ avec $\wt k$ et ainsi remplacer, gr\^ace au th\'eor\`eme de Hukuhara-Turrittin (voir par exemple~\cite{Malgrange91}), $\wt M'$ et $\wt M''$ par leurs mod\`eles \'el\'ementaires respectifs $\oplus_\varphi (E^\varphi\otimes R'_\varphi)$ et $\oplus_\varphi (E^\varphi\otimes R''_\varphi)$.

\begin{lemme}
Si $\varphi,\psi\in\ym\CC[\ym]$ sont distincts, tout accouplement sesquilin\'eaire $\wt k_{\varphi,\psi}=(E^\varphi\otimes R'_\varphi)\otimes_\CC(E^{-\ov\psi}\otimes \ov{R''_{-\psi}})\to\Dbmo_{\wt Y}$ prend ses valeurs dans le sous-faisceau des fonctions \`a d\'ecroissance rapide.
\end{lemme}

\begin{proof}
Puisque $e^{\psi-\ov\psi}$ est un multiplicateur sur $\Dbmo_{\wt Y}$, on peut se ramener au cas o\`u, par exemple, $\psi=0$. Par r\'ecurrence sur le rang de $R'_\varphi$ et $R''_0$, on se ram\`ene au cas de rang $1$, et puisque les fonctions $y^\alpha$ ou $\ovy^\beta$ sont aussi des multiplicateurs, on se ram\`ene au cas o\`u $R'_\varphi$ et $R''_0$ sont \'egaux \`a $\cO_Y$. Alors $\wt u=\wt k(\ephi,1)$ est un germe de distribution mod\'er\'ee sur $\wt Y$ qui satisfait \`a $\ov\partial_y\wt u=0$ et $\partial_y\wt u=\varphi'(y)\wt u$. On en d\'eduit $\wt u_{|Y^*}=e^{\varphi}$. Donc $\wt u$ \`a croissance mod\'er\'ee $\iff$ $\wt u$ \`a d\'ecroissance rapide.
\end{proof}

De la m\^eme mani\`ere (en utilisant la forme normale pour $R'_0,R''_0$), on voit que les termes diagonaux $\wt k_{\varphi,\varphi}(\wt m',\ov{\wt m''})$ se d\'ecomposent en somme, \`a coefficients dans $\cC^\infty_{\wt Y}$, de termes $e^{\varphi-\ov\varphi}y^{\beta'}\ovy^{\beta''}(\log y)^j(\log \ovy)^k$ ($\beta',\beta''\in\CC$, $j,k\in\NN$). On r\'e\'ecrit chacun de ces termes comme une somme, \`a coefficients dans $\cC^\infty_{\wt Y}$, de termes $\my^{2\beta}\Ly^\ell$ ($\beta\in\CC$, $\ell\in\NN$).

Si $m',m''$ sont des sections locales de $M',M''$, on utilise une partition de l'unit\'e sur~$\wt Y$. On obtient pour $e^*k(m',\ov{m''})$ un d\'eveloppement du type $(*)$, \`a coefficients $\wt f_{\varphi,\beta,\ell}$ dans $e_*\cC^\infty_{\wt Y}$, \`a l'ajout pr\`es d'une fonction $C^\infty$ infiniment plate le long de $e^{-1}(0)$: on l'incorpore dans un des coefficients $\wt f_{\varphi,\beta,\ell}$. On note alors $\wt B_\varphi$ l'ensemble d'indices $\beta$ correspondant \`a $\varphi$. Puisque $\my$ est $C^\infty$ sur $\wt Y$, on peut supposer que la diff\'erence de deux \'el\'ements distincts de $\wt B_\varphi$ n'est pas dans $\hZZ$.

Il reste \`a voir que l'on peut r\'e\'ecrire ce d\'eveloppement avec des coefficients $f_{\varphi,\beta,\ell}$ dans $\cC^\infty_Y$. Nous allons utiliser un argument de transformation de Mellin, comme dans le corollaire \ref{cor:asympt}, dont nous n'utiliserons que les notations.

Notons $y=\rho e^{i\theta}$ en coordonn\'ees polaires. Une fonction $\wt f\in e_*\cC^\infty_{\wt Y}$ admet un d\'eveloppement de Taylor $\sum_{m\geq0}\wt f_m(\theta)\rho^m$ o\`u $\wt f_m(\theta)$ est $C^\infty$ sur $S^1$ et se d\'eveloppe en s\'erie de Fourier $\sum_n\wt f_{mn}e^{in\theta}$. Une telle fonction peut s'\'ecrire sous la forme $\sum_{k=-2k_0}^0g_k(y)\my^k$ avec $k_0\in\NN$ et $g_k\in\cC^\infty(Y)$ si et seulement si
\begin{equation}\label{eq:cns}
\wt f_{m,n}\neq0\implique \frac{m\pm n}{2}\geq -k_0.
\end{equation}
En effet, si cette condition est satisfaite, on \'ecrit $\wt f_{m,n}e^{in\theta}\rho^m=\wt f_{m,n}y^{k'}\ov y^{k''}$ avec $k'=(m+n)/2$ et $k''=(m-n)/2$, donc $k',k''\geq -k_0$ et $k'+k''\geq0$. Il existe alors un entier $k$ compris entre $-2k_0$ et $0$ tel que $k'-k/2\in\NN$ et $k''-k/2\in\NN$. La partie du d\'eveloppement de Fourier de $\wt f$ correspondant \`a $k$ fix\'e fournit, par Borel, une fonction $g_k(y)\in\cC^\infty(Y)$. La diff\'erence $\wt f-\sum_{k=-2k_0}^0g_k(y)\my^k$ est une fonction $C^\infty$ sur $\wt Y$, infiniment plate le long de $\my=0$. C'est donc aussi une fonction $C^\infty$ sur~$Y$, infiniment plate \`a l'origine. On l'ajoute \`a $g_0$ pour obtenir la d\'ecomposition voulue de~$\wt f$.

La condition \eqref{eq:cns} peut s'exprimer en terme de transform\'ee de Mellin. On remarque en effet que, pour tous $k',k''\in\hNN$ tels que $k'+k''\in\NN$, la transform\'ee de Mellin $s\mto\langle \wt f,\my^{2s}y^{-k'}\ov y^{-k''}\chi\rangle$, qui est holomorphe pour $\reel(s)\gg0$, s'\'etend en une fonction m\'eromorphe sur $\CC$ avec des p\^oles simples au plus contenus dans $\hZZ$. La condition \eqref{eq:cns} est \'equivalente au fait qu'il existe $k_0\in\NN$ tel que, pour tous $k',k''\in\hNN$ avec $k'+k''\in\NN$, les p\^oles de $s\mto\langle \wt f,\my^{2s}y^{-k'}\ov y^{-k''}\chi\rangle$ soient contenus dans l'intersection des ensembles $k_0-1+k'-\hNN^*$ et $k_0-1+k''-\hNN^*$.

En raisonnant par r\'ecurrence descendante sur $\ell$, c'est-\`a-dire aussi sur l'ordre maximal des p\^oles, on conclut qu'une fonction $\sum_{\ell=0}^L\wt f_\ell\Ly^\ell$ \`a coefficients dans $\cC^\infty(\wt Y)$ peut se r\'e\'ecrire sous la forme $\sum_{-2k_0\leq k\leq0}\sum_{\ell=0}^Lg_{k,\ell}(y)\my^k\Ly^\ell$ avec $g_{k,\ell}\in\cC^\infty(Y)$ si et seulement si la m\^eme condition est satisfaite (et les p\^oles sont d'ordre $\leq L+1$).

Maintenant, si $\wt B_0\subset\CC$ est un ensemble fini tel que deux \'el\'ements distincts ne diff\`erent pas d'un demi-entier, une fonction $\wt f=\sum_{\beta\in\wt B_0}\sum_{\ell=0}^L\wt f_{\beta,\ell}\my^{2\beta}\Ly^\ell$ \`a coefficients dans $\cC^\infty(\wt Y)$ se r\'e\'ecrit $\sum_{\beta\in B_0}\sum_{\ell=0}^Lf_{\beta,\ell}\my^{2\beta}\Ly^\ell$ pour un certain ensemble $B_0$, avec $f_{\beta,\ell}\in\cC^\infty(Y)$, si et seulement si il existe un ensemble fini $A_0\subset\CC$ tel que, pour tous $k',k''\in\hNN$ avec $k'+k''\in\NN$, les p\^oles de $s\mto\langle \wt f,\my^{2s}y^{-k'}\ov y^{-k''}\chi\rangle$ soient contenus dans $(A_0+k'-\NN)\cap(A_0+k''-\NN)$.

Enfin, si $\wt f$ admet un d\'eveloppement du type $(*)$ \`a coefficients dans $\cC^\infty(\wt Y)$, la condition ci-dessus appliqu\'ee \`a $\wt f$ est \'equivalente au fait que $\wt f$ peut se r\'e\'ecrire avec des coefficients $f_{0,\beta,\ell}\in\cC^\infty(Y)$, d'apr\`es un analogue \'evident du lemme \ref{lem:sanspole}.

On applique donc ceci \`a $k(m',\ov{m''})$: on voit que la condition sur la transform\'ee de Mellin est satisfaite en utilisant l'existence d'\'equations fonctionnelles de Bernstein pour $m'$ et $m''$, de la m\^eme mani\`ere que dans le corollaire \ref{cor:asympt} et on obtient ainsi le r\'esultat pour les coefficients avec $\varphi=0$. On obtient le r\'esultat pour les autres coefficients en appliquant le m\^eme raisonnement \`a $E^{-\varphi}\otimes M'$ et $E^{\varphi}\otimes M''$ pour les diff\'erents $\varphi$.
\end{proof}

\section{Application \`a la filtration parabolique canonique}\label{sec:fil}
Soit $M$ un germe de fibr\'e m\'eromorphe \`a connexion sur $X$ et $N=\pi^+M$ son image inverse par une ramification $\pi:Y\to X$. On suppose que $N$ admet un mod\`ele \'el\'ementaire formel: $\wh N\simeq\oplus_\varphi(E^\varphi\otimes R_\varphi)$. Pour tout $\varphi$ et tout $b\in\RR$, on note $\ccP^b(\wh N)=\oplus_{\varphi\in\Phi}(E^\varphi\otimes V^b R_\varphi)$, si $V^\cbbullet$ d\'esigne la filtration de Deligne du fibr\'e m\'eromorphe \`a singularit\'e r\'eguli\`ere $R_\varphi$, d\'efinie comme plus haut par le fait que le r\'esidu de la connexion logarithmique sur $V^b R_\varphi$ a ses valeurs propres de partie r\'eelle dans $[b,b+1[$. Alors $\ccP^b(\wh N)$ est un $\CC\lcr y\rcr$-module libre de type fini tel que $\CC\lcr y\rcr[\ym]\otimes_{\CC\lcr y\rcr}\ccP^b(\wh N)=\wh N$. On a aussi $y\ccP^b(\wh N)=\ccP^{b+1}(\wh N)$. Toute suite exacte $0\to \wh N'\to \wh N\to \wh N''\to0$, o\`u $\wh N,\wh N',\wh N''$ admettent une d\'ecomposition comme ci-dessus, est strictement filtr\'ee relativement \`a la filtration~$\ccP^\cbbullet$.

Le $\CC\{y\}$-module $\ccP^b(N)\defin\ccP^b(\wh N)\cap N$ est libre de type fini, c'est un r\'eseau de~$N$ et toute suite exacte $0\to N'\to N\to N''\to0$, o\`u $N,N',N''$ admettent un mod\`ele \'el\'ementaire formel, est strictement filtr\'ee relativement \`a la filtration~$\ccP^\cbbullet$.

Enfin, en prenant les invariants sous l'action de $\ZZ/q\ZZ$, on obtient la \emph{filtration (parabolique) canonique} de $M$, qui se comporte aussi de mani\`ere stricte dans toute suite exacte (\cf \cite{Malgrange95} pour un cas plus g\'en\'eral de cette construction): $\ccP^b(M)\defin \ccP^{qb}(N)^{\ZZ/q\ZZ}$. Remarquons que, pour tout $k\in\ZZ$, on~a $x^k\ccP^b(M)=\ccP^{b+k}(M)$.

\begin{lemme}
Dans $M$ on~a, pour tout $b\in\RR$, l'\'egalit\'e $\CC\{x\}\langle x\partial_x\rangle\cdot \ccP^b(M)=V^b(M)$.
\end{lemme}

\begin{proof}
On remarque d'abord que la filtration $V^\cbbullet$ se comporte comme $\ccP^\cbbullet$ par rapport \`a la ramification, c'est-\`a-dire que $V^b(M)=V^{qb}(N)^{\ZZ/q\ZZ}$: en effet, $y\partial_y$ op\`ere comme $q(\id\otimes x\partial_x)$ sur $1\otimes M\subset N$. L'inclusion $\subset$ est alors claire. On se ram\`ene \`a montrer l'\'egalit\'e pour $N$, en prenant ensuite les invariants sous $\ZZ/q\ZZ$. Il suffit enfin de montrer l'\'egalit\'e des formalis\'es des deux termes, ce qui donnera aussi leur \'egalit\'e.

D'un c\^ot\'e on~a $V^b(\wh N)=V^b(R_0)\oplus\bigoplus_{\varphi\neq0}(E^\varphi\otimes R_\varphi)$. D'un autre c\^ot\'e, on~a
\[
\CC\{y\}\langle y\partial_y\rangle\cdot (E^\varphi\otimes V^b(R_\varphi))=
\begin{cases}
V^b(R_0)&\text{si }\varphi=0,\\
E^\varphi\otimes R_\varphi&\text{sinon.}
\end{cases}
\]
On en d\'eduit l'assertion.
\end{proof}

\begin{corollaire}\label{cor:extmin}
L'extension minimale $M_{\min}\defin \CC\{x\}\langle \partial_x\rangle\cdot V^{>-1}(M)$ est aussi \'egale \`a $\CC\{x\}\langle \partial_x\rangle\cdot \ccP^{>-1}(M)$.\qed
\end{corollaire}

Soient $M',M''$ deux fibr\'es m\'eromorphes sur $X$ et soit $k:M'\otimes_\CC\ov{M''}\to \Dbmo_{X,0}$ un germe d'accouplement $\cD_{X,0}\otimes_\CC\cD_{\ov X,0}$-lin\'eaire.

\begin{corollaire}\label{cor:restrk}
La restriction de $k$ \`a $\ccP^{>-1}(M')\otimes_\CC\ov{\ccP^{>-1}(M'')}$ prend ses valeurs dans l'espace des (germes de) fonctions $L^1_\loc(\itwopi dx\wedge d\ov x)$.
\end{corollaire}

\begin{proof}
Il est \'equivalent et plus commode de montrer que la restriction de $k$ \`a $\ccP^{>0}(M')\otimes_\CC\ov{\ccP^{>0}(M'')}$ prend ses valeurs dans $L^1_\loc(\itwopi \frac{dx}{x}\wedge \frac{d\ov x}{\ov x})$. Soit $\pi:y\mto x=y^q$ une ramification adapt\'ee \`a $M=\nobreak M'$ et~$M''$ et notons $N=\pi^+M$. On peut d\'efinir de mani\`ere naturelle un accouplement $\pi^+k:N'\otimes_\CC\ov{N''}\to\Dbmo_{Y,0}$ qui est $\cD_{Y,0}\otimes_\CC\cD_{\ov Y,0}$-lin\'eaire et dont la restriction aux \'el\'ements $\ZZ/q\ZZ$-invariants soit $\pi^*k$. Il suffit alors de montrer la proposition pour $\pi^+k$ et $L^1_\loc(\itwopi \frac{dy}{y}\wedge \frac{d\ov y}{\ov y})$. Si $m'\in\ccP^{>0}(N')$ et $m''\in\ccP^{>0}(N'')$ alors, en posant $v=\pi^+k(m',\ov{m''})$, tout $\beta\in B_\varphi(v)$ a une partie r\'eelle $>0$ et le corollaire \ref{cor:asympt} appliqu\'e \`a $v$ montre que, puisque $e^{\varphi-\ov\varphi}$ est localement born\'ee pour tout $\varphi\in\ym\CC[\ym]$, $v$ est dans~$L^1_\loc$.
\end{proof}

\begin{remarque}[Filtration parabolique canonique et cycles proches formels]
Pour tout $b\in\RR$, le quotient $\gr^b_\ccP(N)\defin\ccP^b(N)/\ccP^{>b}(N)$ est isomorphe \`a $\oplus_\varphi\oplus_{\beta\mid\reel\beta=b}\psi_y^\beta(R_\varphi)$; il est de dimension finie. Il est muni d'un endomorphisme semi-simple $\rS=\oplus_\varphi\oplus_{\beta\mid\reel\beta=b}\beta\id$ et d'un endomorphisme nilpotent $\rN$ somme directe des endomorphismes $y\partial_y-\beta$ sur les $\psi_y^\beta(R_\varphi)$. Enfin, $y:\psi_y^\beta(R_\varphi)\to\psi_y^{\beta+1}(R_\varphi)$ est un isomorphisme compatible avec l'action de $\rN$.

On pose alors, pour $\reel \beta=b$,
\[
\wh\psi_y^\beta(N)=\ker[\rS-\beta\id:\gr^b_\ccP(N)\ra\gr^b_\ccP(N)]=\oplus_\varphi\psi_y^\beta(R_\varphi).
\]
Cet espace, muni de $\rN$, est l'\emph{espace des cycles proches formels} de $N$ pour la valeur propre $\beta$ (il est plus gros que $\psi_y^\beta N=\psi_y^\beta R_0$).

L'action de $\ZZ/q\ZZ$ sur $N$ pr\'eserve $\ccP^b(N)$ et induit une action sur $\gr^b_\ccP(N)$. Elle pr\'eserve aussi chaque $\oplus_\varphi\psi_y^\beta(R_\varphi)=\wh\psi_y^\beta(N)$, \ie commute \`a $\rS$, et elle commute aussi \`a $\rN$.

En prenant les invariants sous $\ZZ/q\ZZ$, on trouve
\[
\gr^b_\ccP(M)\defin \ccP^b(M)/\ccP^{>b}(M)=\gr^{qb}_\ccP(N)^{\ZZ/q\ZZ},
\]
et on munit le terme de gauche de l'endomorphisme semi-simple $\rS_M=\frac1q\rS_N$ et de l'endomorphisme nilpotent $\rN_M=\frac1q\rN_N$. On en d\'eduit
\[
(\wh\psi_y^\beta(M),\rN)\defin\ker[\rS-\beta\id:\gr^b_\ccP(M)\ra\gr^b_\ccP(M)]=(\wh\psi_y^{q\beta}(N)^{\ZZ/q\ZZ},\tfrac1q\rN).
\]

Toute suite exacte $0\to M'\to M\to M''\to 0$ donne lieu \`a une suite exacte $0\to \wh\psi_y^\beta M'\to \wh\psi_y^\beta M\to \wh\psi_y^\beta M''\to 0$ compatible \`a $\rN$.

On a $(\wh\psi_y^\beta (M),\rN)=(\wh\psi_y^\beta(\wh M),\rN)$ pour tout $\beta$. Si $M$ est r\'egulier en $0$, on~a $(\wh\psi_y^\beta (M),\rN)=(\psi_y^\beta (M),\rN)$ pour tout $\beta$.
\end{remarque}

\backmatter
\providecommand{\bysame}{\leavevmode ---\ }
\providecommand{\og}{``}
\providecommand{\fg}{''}
\providecommand{\smfandname}{\&}
\providecommand{\smfedsname}{\'eds.}
\providecommand{\smfedname}{\'ed.}
\providecommand{\smfmastersthesisname}{M\'emoire}
\providecommand{\smfphdthesisname}{Th\`ese}

\end{document}